\documentclass[a4paper,10pt]{article}
\def \Z {{\mathbf {Z}}}
\def \R {{\mathbf {R}}}
\def \N {{\mathbf {N}}}

\def \B {{\cal B}}

\title{ Entropy invariants of generic  actions}
\author{V.V. Ryzhikov}
\date{}
\begin{document}
\Large
\maketitle
\begin{abstract} We show  that the typical dynamical system sometimes begins 
to behave like a non-deterministic system with a small classical entropy, and this behavior lasts an extremely long time,
until  the system starts decreasing entropy. Then again it will become almost non-deterministic
for a very very long time,  but with  more smaller  classical entropy. Playing on this fact
 and  considering sigma-compact families of measure-preserving zero-entropy transformations,
 for example, the rectangle exchange transformations, we  choose the Kushnirenko entropy 
so that it is equal to zero for the transformations under consideration, but is infinite for the generic transformation.  
\end{abstract} 

The typical dynamical system is not mixing, but at some points in time it becomes like mixing, 
and such mixing intervals may last as long as the reader desires it and even longer. 
This fact was used in \cite{R1} to refine the theorem on  the generic  actions of  lattices \cite{Ti}.
Now we would like to show that some powers of the generic transformation at very long time
behave as  a non-deterministic system.

 In \cite{CD} it was proved    the non-typicalness of the interval exchange  transformations (IETs). 
In connection with this the following  question arised: 
will the rectangle exchange  transformation be non-typical?
Little  is known about  the mixing properties of such systems, except that they have zero 
entropy\footnote {J.-P. Thouvenot provided the author with an elegant proof of this and more general  fact
about the entropy of   piecewise isometric transformations.}.
We  answer positively the question  by use of special  Kushnirenko's   invariants, which are   
infinite  for  the  generic transformations and  zero for the systems under  consideration. In particulare,
we give the entropy proof of the mentioned result  on IETs.

\section{Generic properties of transformations}
 We fix a standard probability space $(X,\B,\mu)$ and consider the group of its automorphism $Aut$, which is 
 equipped by the Halmos complete metric $\rho$.  When we traditionally  write
 "the  generic (typical) transformation is weakly mixing", 
 we mean that there is   a dense $G_\delta$-set consisting only of weakly mixing transformations. We say also that
the set of all weakly mixing transformations is comeager. 
In fact there are generic properties, there are comeager sets, but there is no   generic transformation 
(such a  transformation that possesses all generic properties). For example, if $T$ is weakly mixing, then 
$T^{m_i}\to_w\Theta$, where 
$\Theta$ stands for the orthoprojection to the constants in $L_2(\mu)$.  But this  $T$ does not possess
 the  following generic property: there is a subsequence $i_k$ such that $T^{m_{i_k}}\to I$. 
As the reader has already noted, we use the same notation for the transformation and the corresponding operator
in $L_2(\mu)$.

\vspace{3mm}
\bf 1.1. To be  weakly mixing and not to be mixing. \rm These properties  both are generic  (P. Halmos, V. Rokhlin).  
To see this  again   let us fix an infinite set $M\subset\N$ and
find a comeager set $Y\subset Aut$ such that for all  $T\in Y$  one has  
$$T^{m_i}\to\Theta,\ \ T^{n_i}\to I$$
for some sequences $\{m_i\},\{n_i\}\subset M$.
By this reason \it
for any transformation $T$ the set  of all transformations spectrally disjoint from $T$ 
  is comeager. \rm

Let us prove a bit more general assertion. An operator function  $Q(T)$ is said  admissible,
 if  $Q(T)=a\Theta+\sum_i a_iT^i$ and   $a, a_i\geq 0$, $a+ \sum_i a_i=1$.

\vspace{3mm}
\bf Theorem 1.1. \it For an  infinite set $M$ and any admissible functions   $Q$,$R$ 
\\ the \ set\ \ $\{T: T^m\to_w Q(T)\ as\    m\in M, m\to\infty\}$ is meager,
\\
the \ set\ \ $Y=\{T: T^m\to_w R(T)\ for \ some \ \ \{m_k\}\subset M\}$ is comeager.
\rm  

\vspace{3mm}
\bf Corollary. \it For any transformation $T$ the set  of all transformations spectrally disjoint from $T$ 
  is comeager.   So the orbit $\{RTR^{-1}: R\in Aut\}$  is meager \cite{Ju} (see also \cite{FW}).\rm 

\vspace{3mm}
Proof.  The first claim follows from the second. To get  the latter  we fix  $\{J_q\}$, a dense 
  family of transformations,  then choose $M'\subset M$ and a transformation $S$ such that $S^m\to_w R(S)$,  $m\in M'$.   
It is not hard to construct such $S$ as a rank-one transformation.
For any $n$ and $q$  we find  $m=m(n,q)$ and  a neighbourhood 
$U(n,q)$ of $J_q^{-1}SJ_q$ such that the unequality 
$$w(T^m,R(T))< \frac 1 n$$ 
holds  for all  $T\in U(n,q)$.  We get a dense $G_\delta$-set
$$\bigcap_n\bigcup_q U(n,q)\subset Y.$$  

There are many  applications of such weak limits, see, for example,  \cite{O1}, \cite{R3}, 
and \cite{St}, where the typicalness of the limits 
$a\Theta+ (1-a)I$  has been proven.  

\vspace{3mm}
\bf 1.2. Asymmetry.  \rm If  $T$ and $T^{-1}$ are not isomorphic, we call $T$  asymmetric. 
The first example of such $T$ has been presented by N.Anzai.  A weakly  mixing  asymmetric  transformation
was appeared in  \cite{O2}, the typicalness of this property  has been  established in  \cite{Ju}.
Later it was appeared  the following more special property (see \cite{R2}).

\vspace{3mm}
\bf  
Theorem 1.2. \it  For some sequences $m_i,n_i$ there is a transformation $T$ with the property:
$$\mu(A\cap T^{m_i}A\cap T^{n_i}A)\to (\mu(A)+2\mu(A)^3)/3,$$
$$\mu(A\cap T^{-m_i}A\cap T^{-n_i}A)\to \mu(A)^2.$$
The set of transformations $S$ with the same property  for    subsequences of the sequence   $i$
is comeager. All such $S$ are  asymmetric.\rm 

\vspace{3mm}
\bf 1.3. Refined typicalness. \rm The work \cite{Ki}  on   the  roots of the generic  transformation 
stimulated the discovery of new  facts in 
 ergodic theory of generic  measure-preserving systems.  
A generic  transformation  has many roots, it is a finite group extension 
\cite{Ag}. At the same time it is a relative weakly mixing extension of some nontrivial 
factor \cite{GW}. Generic transfrmations are    embedded  in a  flow  \cite{LR}, 
even in weakly mixing $\R^n$-actions \cite{Ti}. 
The centralizer of the generic transformation has a reach structure, which   contains a free  action of 
the infinite-dimensional torus \cite{SE}.    

The theory of the generic group actions  is not quite a thing in itself, 
sometimes it turns out to be useful for applications.  
O. Ageev used nontrivial generic arguments to solve the  homogeneous spectrum  problem  in the class of the weakly mixing
   transformations, and later S. Tikhonov did the same in case of  mixing.
\section {Non-typical transformations}
\it Let $ K $ be a compact in $(Aut,\rho)$  family of transformations. Is it true that the orbit of $K$, i.e. 
$$ K^{Aut}=\{STS^{-1}\ : T\in K, S\in Aut\},$$ is meager, 
in other words,  out of it there is a dense $ G_\delta$-set? \rm

In \cite{R3} the following partial result is appeared.

\vspace {3mm}
\bf Theorem 2.1. \it Let $ K \subset Aut $ be compact set in $(Aut,\rho)$ and for some $ r> 0 $ for all
transformations $ T \in K $ and any positive integer $ n $, there exists $ m> n $ such that $ w(T^m, \Theta)> r $,
where  $ w $ is a metric defining the weak operator topology.
Then $K^{Aut}$  is meager. \rm
  
{Proof.}
 Let $ min (T, j) $ be the minimal number among those $ m> j $ for which $ w (T^m, \Theta)> r $.
Since $ K $ is compact, $ min (T, j) $ is a bounded function on $ K $.
We denote by $ M (j) $ the maximum value of $ min (T, j) $ and  consider the sets $ F_j = \{j, j + 1, \dots, M (j) \} $.
It was proved in \cite {R1} that there exists comeager set $ Y $ such that for every $ S \in Y $
there is a sequence $ {j (k)} $ such that
$$ w (S^m, \Theta) \to 0 \ as\  m\in  F_{j (k)},  m\to\infty.  $$
This transformation $ S$ does not belong to $ K $, otherwise 
$$ \exists m \in F_{j (k)} \ \ dist (S^m, \Theta)> r, $$
but this is forbidden for the transformatioms from   $Y$. Since $Y^{Aut}\cap K=\phi$, we get $Y\cap K^{Aut}=\phi$.

The interval exchange transformations (as we fix a number of intervals) satisfy the conditions of Theorem 2.1.
Thus, we have a simpler proof of  the result  \cite {CD} about non-typilcalness of IETs. 
Our proof  in shorter, since we did not care about estimating the numbers $ M (j) $. 
The partial rigidity   of the IETs (see \cite{Ka}) provides the condition $ w(T^m, \Theta)> r$.

The remarks kindly sent by Benjamin Weiss suggested to the author the following version of Theorem 2.1.

\vspace{3mm}
\bf Theorem 2.2.  \it The set of the   transformations that are  
disjoint from all ergodic IETs  is comeager. In short,  $\{IETs\}^\perp$ is comeager.\rm

\vspace{3mm}
Proof.   Let $T$ be generic and $S$ be  an ergodic  IET.  
We can find a sequence $m_k$ such that 
$$ T^{m_k}\to\Theta,\ \  S^{m_k}\to aI+ (1-a) P,$$
where $P$ is some Markov operator commuting with $S$. 
The sequence $m_r$ is said mixing for $T$ and patially rigid for $S$.
We find such partially rigid sequences $m_k$ within the 
above mixing sets $F_{j(k)}$ corresponding to the generic transformations $T$ (see the proof of Theorem 2.1).
The disjointness follows easy from the above conditions.   Let  $J$ be  Markov operator such that $SJ=JT.$
Then we have
    $$S^mJ=JT^m,\ \ S^{m_k}J=JT^{m_k},$$ 
$$aJ+ (1-a) PJ=J\Theta=\Theta.$$
Since $\Theta$ is an  extreme point in the set of all Markov operators that intertwine ergidic $S$ with weakly mixing $T$, 
we get 
$$J=\Theta.$$
Thus, $S$ and $T$ are disjoint in sense of Furstenberg.

However, we cannot use similar arguments for rectangle exchange transformations.
The mixing properties of the latters have not been studied.  For this reason, we will use the following entropy
arguments.

\section{Entropy invariants}
Among the wide variety of the entropy notions   
(see, for example, \cite{Bo},\cite{KT},\cite{Th},\cite{VZ})
we consider  slightly modified invariants of Kushnirenko \cite{Ku}.

 Let $P=\{P_j\}$ be a sequence of finite subsets of a coutable infinite 
group G. We suppose that $|P_j|\to\infty$. For
 a measure-preserving $G$-action  $T=\{T_g\}$   we define
$$h_j(T,\xi)=\frac 1 {|P_j|}  H\left(\bigvee_{p\in P_j}T_p\xi\right),$$
$$h_{P}(T,\xi)={\limsup_j} \ h_j(T,\xi),$$
$$h_{P}(T)=\sup_\xi h_{P}(T,\xi),$$
$$h^{inf}_{P}(T)=\sup_\xi\liminf_j \ h_j(T,\xi),$$
where $\xi$ is a finite measurable partition of $X$.

\bf 3.1. Upper $P$-entropy. \rm  Without getting carried away with an overly general situation, 
let's consider below a rather ascetic case of progressions:
  $G=Z$, $P_j=\{j,2j,\dots, L(j)j\}$, for some sequence $L(j)\to\infty$.

\vspace{3mm}
\bf Theorem 3.1. \it The class $\{S: h_P(S)=\infty\}$ is  comeager.\rm

Proof.    Let $\{J_q\}$, $q\in \N$, be dense in  $Aut$,  and
$T$ be a Bernoulli transformation, so are  $T_q=J_q^{-1}TJ_q$.  The set  $\{T_q\}$ is dense in $Aut$.
 We fix a dense collection of partitions $\xi_i$ (in fact the density  is not nessesary).

 For any $n,q$ there is $j=j(n,q)$  such that for all $i\leq n$  
$$ h_{j}(T_q,\xi_i)=\frac 1 {L_j} H(\bigvee_{n=1}^{L(j)} T^{nj}_q\xi_i) >H(\xi_i)-\frac 1 n. \eqno (n)$$
Indeed, $T_q$ is Bernoulli, we find a partition $\xi$ which is close to $\xi_i$ ($i$ is fixed) and 
for some $M(i,q,n,)$  the partitions $T^{nj}_q\xi$ are independent for all $n$ as $j>M(i,n,q)$.
This implies  $(n)$  for sufficiently large $j$.

There is   a neighbourhood $U(n,q)$ of the transformation $T_q$  such that   all $S\in U(n,q)$ sutisfy 
 the same inequality
$$ h_{j}(S,\xi_i) >H(\xi_i)-\frac 1 n.$$ 
We consider 
$$W=\bigcap_n\bigcup_q U(n,q),$$ 
which is a dense $G_\delta$-set. If $S\in W$, then for any $n$ there is $q(n)$ 
such that the inequality   $$ h_{j(n,q(n))}(S,\xi_i)> H(\xi_i)-\frac 1 n$$ holds  for all $i\leq n$. 
This implies  $h_P(S)=\infty$. Thus, the family  $\{S: h_P(S)=\infty\}$ contains the comeager set $W$.

And what about the  group actions? 
The above proof is valid  for the following  generalization.

\vspace{3mm}
\bf Theorem 3.1.1. \it For a given coutable infinite 
group $G$, let $Q$ be a sequence of the collections  $Q_j=\{q_j(1), q_j(2), \dots, q_j(L(j))\}\subset G$, $L(j)\to\infty$,
  such that for any finite $F\subset G$ for all large $j$ and all $m,n,$  $m\neq n\leq L(j)$ the poduct $q_j(m)^{-1}q_j(n)$
is outside $F$.

Then the set of actions $\{T: h_Q(T)=\infty\}$ is  comeager in the space of all $G$-actions.\rm

\vspace{3mm}
\bf 3.2. Compact sets without of generic entropy properties. \rm Recall that 
$K^{Aut}=\{J^{-1}SJ:  S\in K, J\in Aut\}$ (the orbit of $K$). The class of zero entropy transformations 
 is denoted by   $E_0$.

\vspace{3mm}
\bf Theorem 3.2. \it If $K\subset E_0$ is a compact set  in $(Aut,\rho)$,  then $K^{Aut}$ is meager.\rm

\vspace{3mm}
Remark. To prove that the set $K^{Aut}$ is  meager we need to place it in the comeager  set $E_0$.
 Why an inessential, meager set $Aut\setminus E_0$ significantly interferes with the proof that 
$C^{Aut}$ is  meager for an arbitrary compact set $C$?

\vspace{3mm}
Proof.  We fix a dense collection of partitions $\xi_i$.
If   $h(S)=0$, then for any $i$ 
$$h(S^j,\xi)=\lim_{L\to\infty}\frac 1 {L}  H\left(\bigvee_{p=1}^L S^{jp}\xi_i\right)=0.$$
Given $S\in K$ and $j$,  we find a progression sequence  
$$P(S)=\{P_{j}(S)\},  \ P_j(S)=\{j,2j,\dots, L_S(j)j\}$$  
such that 
$$\frac 1 {|P_j(S)|}  H\left(\bigvee_{p\in P_j(S)}S^p\xi_i\right)<\frac 1 j$$
holds  for all $i<j$. 

 Since  $K$ is a compact set,  using  the progression structure of our sets $P_j(S)$, 
we can find a sequence $L(j)$ such that 
  $L(j)>L_S(j)$ for all 
 $S\in K$. For the corresponding progression sequence $P$ we get   $h_P(S)=0$.

 Since the  $\{T: h_P(T)=\infty\}$ 
 is invariant by the conjugations,  we obtain
$$\{T: h_P(T)=\infty\}\cap K^{Aut}=\phi,$$
thus, $K^{Aut}$ is meager.

\vspace{3mm}
\bf Corollary.  \it  The generic transformaton is not an exchange rectangle transformation.\rm

\vspace{3mm}
Indeed, the set of all exchanges of $n$ rectangles  within a fixed rectangle $X$  is compact.
Each of them has zero entropy, what we can explain by use of Jean-Paul's hint.
 Let $\xi$ be a rectangle partition,
then the sum of the boundary measures of atoms in the partition $\xi^N=\bigvee_{n=0}^{N-1} T^n\xi$  grows linearly.
This implies that the entropy of $\xi^N$ grows very slowly and $H(\xi^N)/N \to 0$. 
It is clear that this method
works in a much more general situation. 
Now we apply Theorem 3.2  and get the above  corollary.  

\vspace{3mm}
\bf Remark. \rm J. Chaika   pointed us to the result   \cite{Bu} on topological zero entropy of
piecewise isometric transformations. In fact we get the non-typicalness of several 
natural classes of dynamical systems: for example, 
all explicite classes of rank one transformations,
 or  $\sigma$-compact classes of special flows over rotations. This shows that the ergodists  prefer to 
  work with non-typical systems, the exceptions are  those cases when they 
write articles about generic transformations, including this note.

\vspace{3mm}
\bf 3.3.  Lower $P$-entropy of the generic transfotmation.\rm  

\vspace{3mm}
\bf Theorem 3.3. \it For some progression sequence $P$  the class $\{S: h_P^{inf}(S)=0\}$   
contains a dense $G_\delta$  set.\rm
 
\vspace{3mm}
Proof.    As in the proof of Theorem 3.1 we consider again the partitions $\xi_i$,  
the family  $\{J_q\}$, the dense set  
$T_q=J_q^{-1}TJ_q$ for some ergodic $T$ of zero entropy.  
 For any $n,q$ there is
$j(n,q)$  such that for all $i\leq n$  
$$ h_{j(n,q)}(T_q,\xi_i)<\frac 1 n,$$
and all $S$ from some neighbourhood $U(n,q)$ of  $T_q$   sutisfy the same inequality. We get that
$\bigcap_n\bigcup_q U(n,q)$  is a dense $G_\delta$-subset of the class  $\{S:h_P^{inf}(S)=0\}$. Q.E.D.

\vspace{3mm}
    J.-P. Thouvenot  drew our attention to the preprint  \cite{Ad},
emphasizing some analogy of the  results on entropy typicalness that were independently obtained 
by T. Adams and the author.

\section{Closing remarks}Of course, many questions arise in connection with topological analogs 
and the actions of various groups, 
but we will pay attention to the following problems.

\bf 4.1. The space of mixing actions. \rm
Bashtanov's result \cite{Ba} together with the  well-known facts  shows that the generic mixing transformation 
has no factor and commutes only with its powers. So for  the space  $Mix$ there are no  analogues  of the 
mentioned sophisticated  results by  Ageev, Glasner-Weiss, King, de la Rue-de Sam Lazaro, and Eremenko-Stepin.

Recall that $Mix$ for $\Z$-actions is equipped by 
 Alpern-Tikhonov's metric $d=\rho+m$,
where  $m(S,T)=\sup_{k\in\Z} w(S^k,T^k)$ for a  fixed metric $w$ which defines  the weak operator topology. 
 The space $Mix$ is complete. Indeed,  the Halmos metric $\rho $ is so, then 
from the conditions  $ \rho(T_i,T)\to 0$, $m(T_i,T_j)\to 0$ as $i,j\to\infty$, and 
$T_i\in Mix$ it follows that  $T\in Mix$ and $ m(T_i,T)\to 0$, thus,  $ d(T_i,T)\to 0$.  

\vspace{3mm}
\bf Theorem 4.1. \it By the same reason the  statements of Theorems 3.1, 3.2, 3.3   are valid 
for  the  space  $Mix$. \rm 
  
\vspace{3mm}
\bf 4.2. Asymmetry in Mix.  \rm  Let $T\in E_0$ be multiple mixing. May it possesses the following asymmetry property?
For some sequences $N_k, m_k, n_k\to\infty$ and any generating partition $\xi$  one has 
$$\frac {H(\xi^{N_k}\bigvee T^{m_k}\xi^{N_k}\bigvee T^{n_k}\xi^{N_k})}{H(\xi^{N_k)})}\to a,$$
and for $b\neq a$ 
$$\frac {H(\xi^{N_k}\bigvee T^{- m_k}\xi^{N_k}\bigvee T^{-n_k}\xi^{N_k})}{H(\xi^{N_k})}\to b,$$
where $\xi^N $  denotes the partition $\bigvee_{n=0}^{N-1} T^n\xi$.
 This property (or another  one of a similar nature)
 seems to be generic in $Mix$. It implies that  $T$ is not isomorphic to $T^{-1}$.

\vspace{3mm}
\bf 4.3. $\{2^n\}$-entropy.  \rm 
Classical  Kushnirenko's  $\{2^n\}$-entropy of the horocycle
flow is finite \cite{Ku}, but it is  infinite for the typical transformation.  This fact contrasts with 
the usual impression that the generic measure-preserving transformation has very weak mixing properties, 
while the horocyclic flow has    Lebesgue spectrum and possesses  multiple mixing  property.

Let us consider the sequence $A=\{A_j\}$,
where  $A_j=\{2^j,2^{j+1},\dots, 2^{j^2}\}$.

\vspace{3mm}
\bf Theorem 4.2. \it The family $\{S: h_A(S)=\infty\}$ contains a dense $G_\delta$  set.

\vspace{3mm}
\bf Corollary. \it Kushnirenko's  $\{2^n\}$-entropy of the generic transformation is infinite.  \rm

\vspace{3mm}
In the proof of Theorem 3.1 we replace  $P_j$  by $A_j$. Q.E.D.

\vspace{3mm}
\bf 4.4. Zero $P$-entropy factor. \rm 
A number of questions arise in connection with our attempt to understand 
the entropy properties of a generic transformation.
Does the generic $ T $ have completely positive $P$-entropy?

In connection with the results of Ageev and Glasner-Weiss on  factors
 of  the  generic  transformation,    it would be more interesting  to have in the typical situation ($h_P(T)=\infty$)
 the  existence of a nontrivial factor  with zero $ h_P $-entropy. 


\vspace{3mm}
\bf 4.5. Attractive transformations. \rm Taking a rare opportunity, the author  draws the readers' attention to 
an unusual problem in the theory of typical transformations. 
The question is inspired by the study of analogs of Gordin's homoclinic group.

Let's call a transformation $T$ \it attractive, \rm
if for some transformation $S\neq T $ and some sequence $m_k$  the condition 
$$  T^{-m_k}ST^{m_k}\to T$$ is satisfied.

\it Is the attractiveness   generic?\rm 
\\ In conclusion, we only recall that the Bernoulli actions are
very attractive.  Therefore, in the key of our note, it would be natural to ask,
does the typical transformation  inherit this property?

 \section{Acknowledgements}   
The author thanks T.Adams, L.Bowen, J.Chaika for remarks and questions.
He is also grateful  to  J.-P.Thouvenot and B.Weiss for useful  discussions.


vryzh@mail.ru
\end{document}